\documentclass[12pt]{article}
\usepackage{amsmath}
\usepackage{amssymb}

\usepackage[cp1251]{inputenc}
\usepackage[russian]{babel}
\newcommand{\il}[2]{\int\limits_{#1}^{#2}}
\newcommand{\ilp}[1]{\int\limits_{#1}^{+\infty}}

\newcommand{\ph}{\phantom{a}}
\newcommand{\phh}{\phantom{aaa}}

\newcommand{\sist}[2]{\left\{
\begin{array}{l}
{#1}\\
\ph\\
{#2}
\end{array}
\right.}

\begin{document}

MSC 47B35

\vskip 15pt

\centerline{\bf Convolution type integral equations in the}
\centerline{\bf conservative case}
\vskip 10pt
\centerline{\bf G. A. Grigorian}

\centerline{\it Institute  of Mathematics of NAS of Armenia}
\centerline{\it E -mail: mathphys2@instmath.sci.am}
\vskip 20 pt

Abstract. In this paper convolution  type integral equations in the conservative case are studied. The conservative case of convolution type of equations relates to the case of non normal type of equations and is that of the corresponding symbols degenerate at some points of the real line, and the classical Furier transformation method meets difficulties with its application to the studying equations. To the study in the conservative case of convolution type equations in this paper we use the special factorization method.

\vskip 10pt


{\bf 1. Introduction}. Let $K_j(t), \ph j=1,2$ be absolutely integrable functions on $\mathbb{R}$, and let $f(t)$ be a function on $\mathbb{R}$. Consider the integral equation of two kernels
$$
\phi(t) = f(t) + \ilp{0}K_1(t - \tau) \phi(\tau) d \tau + \il{-\infty}{0}K_2(t - \tau) \phi(\tau) d \tau \eqno (1.1)
$$
and the paired equation
$$
\sist{\phi(t) = f(t) +\il{-\infty}{\infty}K_1(t -\tau) \phi(\tau) d \tau, \ph t < 0,}{\phi(t) = f(t) +\il{-\infty}{\infty}K_2(t -\tau) \phi(\tau) d \tau, \ph t > 0.}  \eqno (1.2)
$$
These equations have important applications in the physics and mechanics and several works are devoted to it (see [1] and cited works in it). For any $K\in L^1(\mathbb{R})$ we set
$$
\nu_m(K) \equiv \il{-\infty}{\infty} t^m K(t) d t, \phh m= 0, 1, 2, \dots .
$$

If $K_j(t) \ge 0, \ph t \in \mathbb{R}, \ph j=1,2,  \ph \nu_0(K_1) = 1$ or $\nu_0(K_2) = 1$ (or both $\nu_0(K_1) =  \nu_0(K_2) = 1$) then the symbols of the equations (1.1) and (1.2) degenerate (vanish) at $0$. This case relates to the case when the application of the classical Furier method to the equations (1.1) and (1.2) is not applicable directly.
In the particular case when
$$
K_j(t) \ge 0, \ph t \in \mathbb{R}, \ph j=1,2, \ph \nu_0(K_1) =  \nu_0(K_2) = 1, \phh \eqno (1.3)
$$

\noindent
$1^o. \ph \nu_1(K_1) < 0, \ph \nu_1(K_2) > 0$

and $f \in L^1(\mathbb{R})$ a result of existence and asymptotic growth of solutions of Eq. (1.2) is obtained by L. G. Arabadjian by the use of nonlinear factorization equations method of N. B. Yengibarian (see [1, p. 221, Theorem 11.2]).

In this paper we interested in whether the special factorization method can be used to the equations (1.1) and  (1.2) in the case $1^o$ and in the (not yet studied) following cases.

\noindent
$2^o. \ph \nu_1(K_1) = 0, \ph \nu_2(K_1) < \infty, \ph   \nu_1(K_2) > 0$

\noindent
$3^o. \ph \nu_1(K_1) = 0, \ph \nu_2(K_1) < \infty, \ph   \nu_1(K_2) < 0$

\noindent
$4^o. \ph \nu_1(K_1) < 0, \ph  \nu_1(K_2) = 0, \ph \nu_2(K_2) < \infty$

\noindent
$5^o. \ph \nu_1(K_1) > 0, \ph  \nu_1(K_2) = 0, \ph \nu_2(K_2) < \infty$

\noindent
$6^o. \ph \nu_1(K_1) > 0, \ph \nu_1(K_2) < 0$

\noindent
$7^o. \ph \nu_1(K_1) < 0, \ph \nu_1(K_2) < 0$

\noindent
$8^o. \ph \nu_1(K_1) > 0, \ph \nu_1(K_2) > 0$

\noindent
provided the conservativeness conditions (1.3) hold. In this paper we use the special factorization method (introduced and developed in the book [2]) to obtain solvability criteria and asymptotic behavior conditions for solutions of the equations (1.1) and (1.2) for the cases $1^o - 8^o$ with the conservativeness conditions (1.3). Notice that the applicability of N. B. Yengibarian's nonlinear factorization equations method to the cases $2^o - 8^o$ with (1.3) remains open.

\vskip 10pt

{\bf 2. Auxiliary propositions}. Denote by $E$ one of the spaces $L^p(\mathbb{R}) \ph (1 \le p < \infty), \linebreak \mathbb{C}^0 \subset \mathbb{M}$, where $\mathbb{M}$ is the space of measurable essentially bounded on $\mathbb{R}$ functions, $\mathbb{C}^0$ is the space of continuous on $\mathbb{R}$ functions $g(t)$ with $\lim\limits_{t \to \pm\infty} g(t) = 0$. Let $k_j(t) \in L^1(\mathbb{R}), \ph j=1,2$ and let $c_j \ne 0, \ph j=1,2$ be any complex numbers. Consider the convolution operators
$$
(A_j\phi)(t) \equiv c_j \phi(t) - \il{-\infty}{\infty}k_j(t - \tau) \phi(\tau) d \tau, \ph \phi \in E, \ph j=1,2
$$
and the functions
$$
\mathcal{A}_j(\lambda) \equiv c_j - \il{-\infty}{\infty} e^{i\lambda t} k_j(t) d t, \ph -\infty \le \lambda \le \infty, \ph j=1,2.
$$
The functions $\mathcal{A}_j(\lambda), \ph j=1,2$  belong to the extended Wiener algebra $W$ and are called the symbols of the operators $A_j, \ph j=1,2$ respectively.  For any $\alpha \in \mathbb{R}$ we set
$$
(B_\alpha\phi)(t) \equiv \phi(t) - (1 + i\alpha) e^t \il{t}{\infty}e^{-s} \phi(s) d s, \ph
(D_\alpha\phi)(t) \equiv \phi(t) - (1 - i\alpha) e^{-t} \il{-\infty}{t}e^s \phi(s) d s,
$$
$t \in \mathbb{R}, \ph \phi \in E$. Let $\alpha_1, \dots, \alpha_r, \ph \beta_1, \dots, \beta_s$ and $\gamma_1, \dots, \gamma_q$ be any real numbers (not all necessarily different) and let $m_1, \dots, m_r, \ph n_1, \dots, n_s$ and $l_1, \dots, l_q$ be any natural numbers (not all necessarily different). Consider the special convolution operators
$$
A_+ \equiv \prod\limits_{j=1}^r D_{\alpha_j}^{m_j}, \phh  A_- \equiv \prod\limits_{j=1}^s B_{\beta_j}^{j_j}, \phh F \equiv \prod\limits_{j=1}^q D_{\gamma_j}^{l_j}.
$$
The functions
$$
\rho_+(\lambda) \equiv \prod\limits_{j=1}^r \biggl(\frac{\lambda - \alpha_j}{\lambda + i}\biggr)^{m_j}, \ph \rho_-(\lambda) \equiv \prod\limits_{j=1}^s \biggl(\frac{\lambda - \beta_j}{\lambda - i}\biggr)^{n_j} \ph \mbox{and} \ph \rho(\lambda) \equiv \prod\limits_{j=1}^q \biggl(\frac{\lambda - \gamma_j}{\lambda + i}\biggr)^{l_j},
$$
$-\infty \le \lambda \le \infty$ are the symbols of the operators $A_+, \ph A_-$ and $F$ respectively (see [2], p.~ 201). Let $P$ be an orthogonal projector, acting on $E$ by the rule
$$
(P\phi)(t) \equiv \sist{\phi(t), \ph t > 0,}{0, \ph t < 0}
$$
and let $Q\equiv I - P$ be another orthogonal projector on $E$, where $I$ is the identity operator on $E$. We set
$$
B \equiv P A_- + Q A_+, \phh D \equiv A_+ P + A_- Q
$$
and denote by $D^{(-1)}$ a linear operator with the properties
$$
(D^{(-1)} D) x = x, \ph x \in E, \ph (D D^{(-1)}) x = x , \ph D^{(-1)} x \in E.
$$
Consider the sets
$$
\widetilde{E}(\rho_+,\rho_-) \equiv D^{(-1)} E, \phh \overline{E}(\rho_+,\rho_-) \equiv im B, \phh \overline{E}(\rho) \equiv im F.
$$
These sets are Banach spaces endowed by some norms, generated by the norm of $E$ and by the operators $A_\pm$ and $F$ as follows
\pagebreak
$$
||f||_{\widetilde{E}(\rho_+,\rho_-)}\equiv ||D f||_E, \phh ||\psi||_{\overline{E}(\rho_+,\rho_-)} \equiv \inf\limits_{B\psi = \phi}||\psi||_E \ph (\phi \in im B),
$$
$$
||\psi||_{\overline{E}(\rho)} \equiv \inf\limits_{F\psi = \phi}||\psi||_E \ph (\phi \in im F)
$$
respectively (see [2], pp. 132, 133, 135). Assume
$$
\mathcal{A}(\lambda) = \rho_+(\lambda) \rho(\lambda) \mathcal{A}_0(\lambda), \ph \mathcal{B}(\lambda) = \rho_-(\lambda) \rho(\lambda) \mathcal{B}_0(\lambda), \ph -\infty \le \lambda \le \infty,
$$
where $\mathcal{A}_0(\lambda), \ph  \mathcal{B}_0(\lambda) \in W$. We set $A = A_1 P + A_2 Q$.

{\bf Theorem 2.1.([2, p. 203, Theorem 4.2].)} {\it In order that $A : \widetilde{E}(\rho_+,\rho_-) \rightarrow \overline{E}(\rho)$ was $\Phi_+$ or $\Phi_-$ operator it is necessary and sufficient that $\mathcal{A}_0(\lambda) \ne 0, \ph \mathcal{B}_0(\lambda) \ne 0, \ph -\infty \le \lambda \le~ \infty.$ If this condition is satisfied, then
$$
dim Ker A = \max \{\delta -\kappa,0\}, \phh \phh dim coKer = \max\{\kappa - \delta,0\}
$$
where $\kappa \equiv \frac{1}{2\pi}[arg (\mathcal{A}_0(\lambda)/\mathcal{B}_0(\lambda))]_{-\infty}^{\infty},$
$$
\delta \equiv \sist{0, \ph if \ph E = L^p \ph (1\le p < \infty) \ph or \ph E = \mathbb{C}^0,}{q \ph if \ph E =M.}
$$

\phantom{aaaaaaaaaaaaaaaaaaaaaaaaaaaaaaaaaaaaaaaaaaaaaaaaaaaaaaaaaaaaaaaaaaaa} $\blacksquare$

A similar result can be formulated for the operator $A' \equiv P A_1 + Q A_2$, which is

{\bf Theorem 2.2.([2, p. 199, Theorem 3.4].)} {\it In order that $A' :E \rightarrow \overline{E}(\rho_+,\rho_-)$ was $\Phi_+$ or $\Phi_-$ operator it is necessary and sufficient that $\mathcal{A}_0(\lambda) \ne 0, \ph \mathcal{B}_0(\lambda) \ne 0, \ph -\infty \le \lambda \le~ \infty.$ If this condition is satisfied, then
$$
dim Ker A' = \max \{\kappa +\sigma,0\}, \phh \phh dim coKer = \max\{ - \kappa - \sigma,0\}
$$
where $\kappa \equiv \frac{1}{2\pi}[arg (\mathcal{A}_0(\lambda)/\mathcal{B}_0(\lambda))]_{-\infty}^{\infty},$
$$
\sigma \equiv \sist{0, \ph if \ph E = L^p \ph (1\le p < \infty) \ph or \ph E = \mathbb{C}^0,}{r + s \ph if \ph E =M.}
$$

\phantom{aaaaaaaaaaaaaaaaaaaaaaaaaaaaaaaaaaaaaaaaaaaaaaaaaaaaaaaaaaaaaaaaaaaa} $\blacksquare$

Let $K(t) \ge 0, \ph t \in \mathbb{R}, \ph K \in L^1(\mathbb{R}).$ Consider the function
$$
a(\lambda) \equiv 1 - \il{-\infty}{\infty}e^{i \lambda t} K(t) d t, \phh -\infty \le \lambda \le \infty.
$$

{\bf Theorem 2.3 ([3, Theorem 2.I])} {\it If $\nu_0(K) =1$ and $\nu_1(K) \ne 0$, then
$$
a(\lambda) = \frac{\lambda}{\lambda + i} b(\lambda), \ph b(\lambda) \ne 0, \ph -\infty \le \lambda \le \infty, \ph b(\lambda) \in W
$$
and
$$
ind \hskip 3pt b(\lambda) = \sist{0, \ph if \ph \nu_1(K) > 0,}{-1, \ph if \ph \nu_1(K) < 0.}
$$

\phantom{aaaaaaaaaaaaaaaaaaaaaaaaaaaaaaaaaaaaaaaaaaaaaaaaaaaaaaaaaaaaaaaaaaaa} $\blacksquare$

{\bf Theorem 2.4 ([3, Theorem 2.II])} \it If $\nu_0(K) =1, \ph \nu_1(K) = 0, \ph \nu_2(K) \le \infty$, then
$$
a(\lambda) = \biggl(\frac{\lambda}{\lambda + i}\biggr)^2 b(\lambda), \ph b(\lambda) \ne 0,  \ph -\infty \le \lambda \le \infty, \ph   b(\lambda) \in W, \ph ind \hskip 3pt b(\lambda) = - 1.
$$

\phantom{aaaaaaaaaaaaaaaaaaaaaaaaaaaaaaaaaaaaaaaaaaaaaaaaaaaaaaaaaaaaaaaaaaaa} $\blacksquare$

{\bf 3. Main results.} Hereafter we will assume that
$$
K_j(t) \ge 0, \ph t \in \mathbb{R}, \ph j=1.2, \phh \nu_0(K_1) = \nu_0(K_2) = 1.
$$
With these conditions the equations (1.1) and (1.2) become equations of non normal type, since their symbols vanish at $0$. In this section we use the theorems of the previous section to prove solvability and asymptotic behavior criteria
for the equations (1.1) and (1.2).

{\bf Theorem 3.1.} {\it The following assertions are valid.

\noindent
$I$ If $\nu_1(K_1) \ne 0, \ph \nu_1(K_2) \ne 0$, then for every $f\in E$ Eq. (1.1) has a solution in $\widetilde{E}\bigl(\frac{\lambda}{\lambda + i}, \frac{\lambda}{\lambda - i}\bigr)$. The corresponding homogeneous equation has:

\noindent
$I_1$ the unique (up to an arbitrary multiplier) nontrivial solution  in $\widetilde{E}\bigl(\frac{\lambda}{\lambda + i}, \frac{\lambda}{\lambda - i}\bigr)$ provided $\nu_1(K_1) > 0, \ph \nu_1(K_2) > 0$ or  $\nu_1(K_1) < 0, \ph \nu_1(K_2) < 0$,

\noindent
$I_2$ only the trivial solution  in $\widetilde{E}\bigl(\frac{\lambda}{\lambda + i}, \frac{\lambda}{\lambda - i}\bigr)$ if  $\nu_1(K_1) > 0, \ph \nu_1(K_2) < 0$,

\noindent
$I_3$ two linearly independent solutions in in $\widetilde{E}\bigl(\frac{\lambda}{\lambda + i}, \frac{\lambda}{\lambda - i}\bigr)$ if  $\nu_1(K_1) < 0, \ph \nu_1(K_2) > 0$.

\noindent
$II$ If $\nu_1(K_1) = 0, \ph \nu_2(K_1) < \infty, \ph \nu_1(K_2) \ne 0$, then for every $f \in E$ Eq. (1.1) has a solution in  $\widetilde{E}\bigl(\bigl[\frac{\lambda}{\lambda + i}\bigr]^2, \frac{\lambda}{\lambda - i}\bigr)$. The corresponding homogeneous equation has:

\noindent
$II_1$ the unique (up to an arbitrary multiplier) nontrivial solution  in $\widetilde{E}\bigl(\bigl[\frac{\lambda}{\lambda + i}\bigr]^2, \frac{\lambda}{\lambda - i}\bigr)$ if $\nu_1(K_2) < 0,$

\noindent
$II_2$ two linearly independent solutions in in $\widetilde{E}\bigl(\bigl[\frac{\lambda}{\lambda + i}\bigr]^2, \frac{\lambda}{\lambda - i}\bigr)$ if $\nu_1(K_2) > 0.$

\noindent
$III$ If $\nu_1(K_1) \ne 0, \ph \nu_1(K_2) = 0, \ph \nu_2(K_2) < \infty,$ then for every $f \in E$ Eq. (1.1) has a solution in  $\widetilde{E}\bigl(\frac{\lambda}{\lambda + i}, \bigl[\frac{\lambda}{\lambda - i}\bigr]^2\bigr)$. The corresponding homogeneous equation has:

\noindent
$III_1$ the unique (up to an arbitrary multiplier) nontrivial solution  in $\widetilde{E}\bigl(\frac{\lambda}{\lambda + i}, \bigl[\frac{\lambda}{\lambda - i}\bigr]^2\bigr)$, if $\nu_1(K_1) > 0$,

\noindent
$III_2$ two linearly independent solutions  in $\widetilde{E}\bigl(\frac{\lambda}{\lambda + i}, \bigl[\frac{\lambda}{\lambda - i}\bigr]^2\bigr)$, if $\nu_1(K_1) < 0$.

\noindent
$IV$ If $\nu_1(K_1) = \nu_1(K_2) = 0, \ph \nu_2(K_j) < \infty, \ph j=1,2,$ then for every $f\in E$ Eq. (1.1) has a solution  $\phi \in \widetilde{E}\bigl(\bigl[\frac{\lambda}{\lambda + i}\bigr]^2, \bigl[\frac{\lambda}{\lambda - i}\bigr]^2\bigr).$
 If $f \in \overline{E}\bigl(\frac{\lambda}{\lambda + i}\bigr)$, then $\phi \in  \widetilde{E}\bigl(\frac{\lambda}{\lambda + i}, \frac{\lambda}{\lambda - i}\bigr)$. The corresponding homogeneous equation has two linearly independent solutions $\phi_1, \ph \phi_2 \in  \widetilde{E}\bigl(\bigl[\frac{\lambda}{\lambda + i}\bigr]^2, \bigl[\frac{\lambda}{\lambda - i}\bigr]^2\bigr).$
If $E = M$, then $\phi_1, \ph \phi_2 \in \widetilde{E}\bigl(\frac{\lambda}{\lambda + i}, \frac{\lambda}{\lambda - i}\bigr)$ and if $E = L^p \ph (1\le p < \infty)$ or $E = \mathbb{C}^0$ then the corresponding homogeneous equation has the unique (up to an arbitrary multiplier) nontrivial solution in $\widetilde{E}\bigl(\frac{\lambda}{\lambda + i}, \frac{\lambda}{\lambda - i}\bigr)$.
}

Proof.  Let $a_1(t)$ and $a_2(t)$ be symbols of the convolution operators
$$
(\mathbf{a}_1\phi)(t)\equiv \phi(t) - \il{-\infty}{\infty}K_1(t - \tau)\phi(\tau) d \tau, \ph  (\mathbf{a}_2\phi)(t)\equiv \phi(t) - \il{-\infty}{\infty}K_2(t-\tau)\phi(\tau) d \tau, \ph   t \in E,
$$
respectively, i. e. $a_j(\lambda) = 1 - \il{-\infty}{\infty}e^{i\lambda t} K_j(t) d t, \ph -\infty \le \lambda \le \infty, \ph j=1,2.$ Let us prove I. Since $\nu_1(K_1) \ne 0, \ph \nu_1(K_2) \ne 0$ by Theorem 2.3 we have $a_1(\lambda) = \frac{\lambda}{\lambda + i} b_1(\lambda), \ph a_2(\lambda) = \frac{\lambda}{\lambda - i} b_2(\lambda), \ph b_j(\lambda) \ne 0, \ph -\infty \le \lambda \le \infty, \ph b_j(\lambda) \in W, \ph j=1,2$ and
$$
ind \hskip 3pt b_1(\lambda) = \sist{0, \ph if \ph \nu_1(K_1) > 0,}{-1 \ph if \ph \nu_1(K_1) < 0,} \ph ind \hskip 3pt b_2(\lambda) = \sist{1, \ph if \ph \nu_1(K_1) > 0,}{0 \ph if \ph \nu_1(K_1) < 0}
$$
(since $a_2(\lambda) = \frac{\lambda}{\lambda + i}\bigl(\frac{\lambda + i}{\lambda - i}\bigr), \ph -\infty \le \lambda \le \infty$ and $ind \hskip 3pt \bigl[\bigl(\frac{\lambda + i}{\lambda - i}\bigr)\bigr] = \sist{0, \ph if \ph \nu_1(K_1) > 0,}{-1 \ph if \ph \nu_1(K_1) < 0} = -1 + ind \hskip 3pt b_2(\lambda)$). Hence,
$$
ind \hskip 3pt \frac{b_1(t)}{b_2(t)} =  \left\{\begin{array}{l}
-1, \ph if \ph \nu_1(K_1) > 0, \ph \nu_1(K_2) > 0,\\
-2, \ph if \ph \nu_1(K_1) < 0, \ph \nu_1(K_2) > 0,\\
0, \ph if \ph \nu_1(K_1) > 0, \ph \nu_1(K_2) < 0,\\
-1, \ph if \ph \nu_1(K_1) < 0, \ph \nu_1(K_2) < 0.
\end{array}
\right.
$$
Therefore (if we take $\rho(\lambda) \equiv 1, \ph \rho_\pm(\lambda) = \frac{\lambda}{\lambda \pm i}$) by Theorem 2.1 the assertion I is fulfilled.

Let us prove II. Since $\nu_1(K_1) = 0, \ph \nu_2(K_1) < \infty,$ and $\nu_1(K_2) \ne 0$ by Theorem 2.4 we have $a_1(\lambda) = \bigl[\frac{\lambda}{\lambda + i}\bigr]^2 b_1(t), \ph b_1(\lambda) \ne 0, \ph -\infty \le \lambda \le \infty, \ph b_1(\lambda) \in W, \ph ind \hskip 3pt b_1(\lambda) = -1,$ and by Theorem 2.3 $a_2(\lambda) = \frac{\lambda}{\lambda - i} b_2(\lambda), \ph b_2(\lambda) \ne 0, \ph -\infty \le \lambda \le \infty, \ph b_2(\lambda) \in W$ and
$$
ind \hskip 3pt b_2(\lambda) = \sist{1, \ph if \ph \nu_1(K_2) > 0,}{0, \ph if \ph \nu_1(K_2) < 0}.
$$
Hence,
$$
ind \hskip 3pt \frac{b_1(\lambda)}{b_2(\lambda)} = \sist{-2, \ph if \ph \nu_1(K_2) > 0,}{-1, \ph if \ph \nu_1(K_2) < 0}.
$$
By Theorem 2.1 from here it follows  II (for this case we take $\rho(\lambda) \equiv 1, \ph \rho_+(\lambda) = \bigl[\frac{\lambda}{\lambda +1}\bigr]^2, \ph \rho_-(\lambda) = \frac{\lambda}{\lambda - i}$).

Let us prove III. Since $\nu_1(K_1) \ne 0, \ph \nu_1(K_2) = 0, \ph \nu_2(K_2) < \infty$ by Theorem 2.3 we have
$a_1(\lambda) = \frac{\lambda}{\lambda + i} b_1(t), \ph b_1(\lambda) \ne 0, \ph -\infty \le \lambda \le \infty, \ph b_1(\lambda) \in W,$
$$
ind \hskip 3pt b_1(\lambda) = \sist{0, \ph if \ph \nu_1(K_1) > 0,}{-1, \ph if \ph \nu_1(K_1) < 0}.
$$

and by Theorem 2.4 $a_2(\lambda) = \bigl[\frac{\lambda}{\lambda - i}\bigr]^2 b_2(\lambda), \ph b_2(\lambda) \ne 0, \ph -\infty \le \lambda \le \infty, \ph b_2(\lambda) \in W$ and $ind \hskip 3pt b_2(\lambda) = 1$
Hence,
$$
ind \hskip 3pt \frac{b_1(\lambda)}{b_2(\lambda)} = \sist{-1, \ph if \ph \nu_1(K_2) > 0,}{-2, \ph if \ph \nu_1(K_2) < 0.}
$$
By Theorem 2.1 from here it follows III (here we take $\rho(\lambda) \equiv 1, \ph \rho_+(\lambda) = \frac{\lambda}{\lambda + i}, \ph \rho_-(\lambda) = \bigl[\frac{\lambda}{\lambda - i}\bigr]^2$). It remains to prove IV. Since $\nu_1(K_1) = \nu_1(K_2) = 0, \ph \nu_2(K_j) < \infty, \ph j=1.2$ by Theorem 2.4 we have
 $a_1(\lambda) = \bigl[\frac{\lambda}{\lambda - i}\bigr]^2 b_1(\lambda), \ph b_1(\lambda) \ne 0, \ph -\infty \le \lambda \le \infty, \ph b_2(\lambda) \in W$,
$a_2(\lambda) = \bigl[\frac{\lambda}{\lambda - i}\bigr]^2 b_2(\lambda), \ph b_2(\lambda) \ne 0, \ph -\infty \le \lambda \le \infty, \ph b_2(\lambda) \in W$ and $ind \hskip 3pt b_1(\lambda) = -1, \ph ind \hskip 3pt b_1(\lambda) = 1$. Then
$$
ind \hskip 3pt \frac{b_1(\lambda)}{b_2(\lambda)} = -2. \eqno (3.1)
$$
Therefore if we take $\rho(\lambda) \equiv 1, \ph \rho_\pm(\lambda) = \bigl[\frac{\lambda}{\lambda \pm i}\bigr]^2$, then by Theorem 2.1 we obtain that Eq. (1.1) has a solution $\phi \in \widetilde{E}\bigl(\bigl[\frac{\lambda}{\lambda + i}\bigr]^2, \bigl[\frac{\lambda}{\lambda - i}\bigr]^2)$. If we take $\rho(\lambda) = \frac{\lambda}{\lambda}{\lambda + i}, \ph \rho_\pm(\lambda) = \frac{\lambda}{\lambda\pm I}$, the by Theorem 2.1 from (3.1) we obtain that $\phi\in \widetilde{E}\bigl(\frac{\lambda}{\lambda + i}, \frac{\lambda}{\lambda - i}\bigr)$ for $f \in \overline{E}\bigl(\frac{\lambda}{\lambda + i}\bigr)$. Moreover,
$$
Ind \hskip 3pt A = \sist{2, \ph for \ph E = M,}{1, \ph for \ph E = L^p(1\le p < \infty) \ph or \ph E = \mathbb{C}^0.}
$$
Therefore, if $E = M$, then $\phi_1, \ph \phi_2 \in  \widetilde{E}\bigl(\frac{\lambda}{\lambda + i}, \frac{\lambda}{\lambda - i}\bigr),$ and if $E = L^p(1\le p < \infty) \ph or \ph E = \mathbb{C}^0$, then corresponding to (1.1) homogeneous equation has the unique (ut to an arbitrary multiplier) nontrivial solution in $\widetilde{E}\bigl(\frac{\lambda}{\lambda + i}, \frac{\lambda}{\lambda - i}\bigr).$ The assertion IV is proved. Thus the proof of the theorem is completed.

Using Theorem 2.2 instead of Theorem 2.1 by analogy with the proof of Theorem 3.1 it can be proved.

{\bf Theorem 3.2.} {\it The following assertions are valid.

\noindent
$I^*$ If $\nu_1(K_1) \ne 0, \ph \nu_1(K_2) \ne 0$, then Eq. (1.2) has a solution in $M$ for $f\in \overline{M}\bigl(\frac{\lambda}{\lambda + i}, \frac{\lambda}{\lambda - i}\bigr).$ The corresponding homogeneous equation has:

\noindent
$I^*_1$ only a trivial solution in $M$ if $\nu_1(k_1) < 0, \ph \nu_1(K_2) > 0$,

\noindent
$I^*_2$ the unique (up to an arbitrary multiplier) nontrivial solution in $M$ if $\nu_1(K_1) > 0, \ph \nu_1(K_2) > 0$ or $\nu_1(k_1) < 0, \ph \nu_1(K_2) < 0$,

\noindent
$I^*_3$ two linearly independent solutions in $M$ if $\nu_1(K_1) > 0, \ph \nu_1(K_2) < 0$.

\noindent
$II^*$ If $\nu_1(K_1) = 0, \ph \nu_2(K_1) < \infty, \ph \nu_1(K_2) \ne 0$, then Eq. (1.2) has a solution in $M$ for $f \in \overline{M}\bigl(\bigl[\frac{\lambda}{\lambda + i}\bigr]^2, \frac{\lambda}{\lambda - i}\bigr).$ The corresponding homogeneous equation has:

\noindent
$II^*_1$ the unique (up to an arbitrary multiplier) nontrivial solution in $M$ if $\nu_1(K_2) > 0$,

\noindent
$II^*_2$ two linearly independent solutions in $M$ if $\nu_1(K_2) < 0$.

\noindent
$III^*$ If $\nu_1(K_1) \ne 0, \ph \nu_1(K_2) = 0, \ph \nu_2(K_2) < \infty$, then Eq. (1.2) has a solution in $M$ for $f \in  \overline{M}\bigl(\frac{\lambda}{\lambda + i}, \bigl[\frac{\lambda}{\lambda - i}\bigr]^2\bigr).$ The corresponding homogeneous equation has:

\noindent
$III^*_1$ the unique (up to an arbitrary multiplier) nontrivial solution in $M$ if $\nu_1(K_1) < 0$,

\noindent
$III^*_2$ two linearly independent solutions in $M$ if $\nu_1(K_1) > 0$.

\noindent
$IV^*$ If $\nu_1(K_1) = \nu_1(K_2) = 0, \ph \nu_2(K_j) < \infty, \ph j=1,2,$ then Eq (1.2) has a solution in $M$ for $f \in  \overline{M}\bigl(\bigl[\frac{\lambda}{\lambda + i}\bigr]^2, \bigl[\frac{\lambda}{\lambda - i}\bigr]^2\bigr) \cup \overline{M}\bigl(\frac{\lambda}{\lambda + i}, \frac{\lambda}{\lambda - i}\bigr).$ The corresponding homogeneous equation has  the unique (up to an arbitrary multiplier) nontrivial solution in $M$.
}

\phantom{aaaaaaaaaaaaaaaaaaaaaaaaaaaaaaaaaaaaaaaaaaaaaaaaaaaaaaaaaaaaaaaaaaaaaa} $\blacksquare$.

\hskip 20pt

\centerline{\bf References}

\hskip 20pt

\noindent
1. L. G. Arabadxhan, N. B. Engibarian, Equations in convolution and nonlinear functional \linebreak \phantom{a} equations (in Russian). Math. Analiz, Itogi Nauki i Tekhnjki, VINITI, Moskow,\linebreak \phantom{a} vol. 22, 1984, pp. 175-244.

\noindent
2. S. Presdorf, Some Classes of Singular Equations, Mir, 1979.

\noindent
3. G. A. Grigorian, Solvability of a class of Wiener-Hopf integral equations. Izv. Nats. \linebreak \phantom{a} Akad. Nauk. Armenii. Matematika, vol 21, No. 2, 1996, pp. 21-32.

\end{document}